\documentclass[a4paper,8pt]{article}
\usepackage{amsmath}
\usepackage{amsthm,amsfonts,amssymb}
\usepackage[dvips]{graphicx}
\usepackage[english]{babel}
\usepackage{color}
\usepackage{graphics}
\usepackage{graphpap}

\providecommand{\U}[1]{\protect\rule{.1in}{.1in}}

\providecommand{\U}[1]{\protect\rule{.1in}{.1in}}
\newtheorem{theorem}{Theorem}[section]

\newtheorem{corollary}[theorem]{Corollary}
\newtheorem{definition}[theorem]{Definition}

\newtheorem{proposition}[theorem]{Proposition}

\makeatletter
\def\blfootnote{\xdef\@thefnmark{}\@footnotetext}
\makeatother

\begin{document}

\title{Transient analysis of the subordinated chain of a state dependent pure birth process}
\author{Andrea Monsellato}
\maketitle

\begin{abstract}
\noindent
Consider a pure birth process with intensities $\lambda_{k}=\frac{1}{1+k}$, with $k=0,1,2,...$, we show that the subordinated chain is assimilable to a Bernoullian scheme with dependent successes probabilities, also we show a direct link with a degenerate Polya urn replacement scheme \cite{Janson}.
We compute explicitly transition probabilities, by generating function method, of the subordinated chain and give some interesting bounds using the \emph{centering sequence} approach proposed by MacDiarmid \cite{McDiarmidcentering}.
\end{abstract}


\section{Introduction}

\noindent
A pure birth process with intensities $\lambda_0,\lambda_1,...\geq 0$ is a process $\{X(t):t\geq 0\}$ taking value in
$S=\{0,1,2,...\}$ such that:\\

\noindent
$X(0)\geq 0$; if $s<t$ then $X(s)\leq X(t)$,\\

\noindent
$P(X(t+h)=k+m|X(t)=k)=\left\{ \begin{matrix} \lambda_mh+o(h), &\mbox{if}\quad m=1\\ o(h),  &\mbox{if}\quad m>1 \\1-\lambda_mh+o(h), &\mbox{if} \quad m=0\end{matrix}\right.$\\

\noindent
if $s<t$ then, conditional on the value of $X(s)$, the increment $X(t)-X(s)$ is independent of all arrivals prior to $s$.\\

\noindent
If $\lambda_n=\lambda$ for all $n$ a birth process with intensity $\lambda_0,\lambda_1,...$ is called \emph{Poisson process}.\\

\noindent
Now consider a pure birth process with intensities $\lambda_k=\frac{\lambda}{1+k}$ with $k\geq 0$ where $\lambda>0$.\\
\noindent
Let $p(k,t)=P(X(t)=k)$ where $\{X(t):t\geq 0\}$ then

\begin{align}\label{equ:bc}
\left\{
\begin{array}{ll}
p'(0,t)=\lambda p(0,t)\\
p'(k,t)= \frac{\lambda}{k}p(k-1,t)-\frac{\lambda}{1+k}p(k,t),\quad if\quad k\neq 0
\end{array}
\right.
\end{align}

\noindent
In the following, without loss of generality, we choose $\lambda=1$.\\

\noindent
It is known that this problem has a unique positive solution and because $\sum_{i=0}^{+\infty}\frac{1}{\lambda_i}=+\infty$ the solution is a proper probability distribution see \cite{Fellervol1}.\\

\newpage
\noindent
\begin{proposition}\label{prop:bpcontsol}

\noindent
Let

\begin{align}\label{solbcont}
\left\{
\begin{array}{ll}
p(0,t)=e^{-t},\\
p(k,t)=\frac{1}{k!}\sum_{j=1}^{k+1}(-1)^{k+1-j}j^k \binom{k+1}{j} e^{-\frac{t}{j}}, \quad if\quad k\geq 1
\end{array}
\right.
\end{align}

\noindent
then (\ref{solbcont}) is the solution of (\ref{equ:bc}).\\

\end{proposition}

\noindent
\textbf{Proof}\\

\noindent
Obviously $p(0,t)=e^{-t}$. Now taking Laplace transform

\begin{equation}\label{lapltransf}
\hat{p}(k,\theta):=\int_{0}^{+\infty}e^{-\theta t}p(k,t)dt
\end{equation}

\noindent
from (\ref{equ:bc}) we have

\begin{align*}
\left\{
\begin{array}{ll}
\hat{p}(0,\theta)=\frac{1}{\theta + 1}\\
\hat{p}(k,\theta)=\frac{1}{k(\theta+\frac{1}{k+1})}\hat{p}(k-1,\theta),\quad k\in [1,n]
\end{array}
\right.
\end{align*}

\noindent
so that

\begin{equation*}
\hat{p}(k,\theta)=\frac{1}{k!}\left[\frac{A_1}{\theta+1}+\frac{A_2}{\theta+
\frac{1}{2}}+\cdot\cdot\cdot+\frac{A_{k+1}}{\theta+\frac{1}{k+1}}\right]
\end{equation*}

\noindent
where $A_i$, $i=1,..,k+1$, satisfy the following equations

\begin{equation*}
A_1\prod_{i=1;i\neq1}^{k+1}\left(\theta+\frac{1}{i}\right)+A_2\prod_{i=1;i\neq 2}^{k+1}\left(\theta+\frac{1}{i}\right)
+A_{k+1}\prod_{i=1;i\neq k+1}^{k+1}\left(\theta+\frac{1}{i}\right)=1
\end{equation*}

\noindent
Choosing $\theta=-\frac{1}{j}$, $j=1,...,k+1$, we obtain:

\begin{align}\label{equ:ajcoeff}
\left\{
\begin{array}{ll}
A_l\prod_{i=1;i\neq j}^{k+1}\left(-\frac{1}{j}+\frac{1}{i}\right)=0,\quad \forall j\neq l\\
\\
A_j=\frac{1}{\prod_{i=1;i\neq j}^{k+1}(-\frac{1}{j}+\frac{1}{i})}
\end{array}
\right.
\end{align}

\noindent
Anti-transforming then for $k\geq 1$

\begin{align*}
p(k,t)=\frac{1}{k!}\sum_{j=1}^{k+1}A_je^{-\frac{t}{j}}=\frac{1}{k!}\sum_{j=1}^{k+1}j^k e^{-\frac{t}{j}}\prod_{i=1;i\neq j}^{k+1}\frac{i}{j-i}
\end{align*}

\noindent
From (\ref{equ:ajcoeff}) we know that

\begin{equation*}
A_j=j^k\prod_{i=1;i\neq j}^{k+1}\frac{i}{j-i}=j^k(-1)^{k+1-j}\prod_{i=1;i\neq j}^{k+1}\frac{i}{|j-i|}
\end{equation*}

\noindent
Let

\begin{equation}\label{Qj}
Q_j:=\prod_{i=1;i\neq j}^{k+1}\frac{i}{|j-i|}
\end{equation}

\noindent
Taking the logarithm of (\ref{Qj}) after some algebraic calculations we obtain

\begin{align*}
Q_j=\frac{(k+1)!}{j!(k+1-j)!}=\binom{k+1}{j}
\end{align*}

\noindent
then the thesis follows.

\begin{flushright}
$\Box$
\end{flushright}

\section{Subordinated chain}

\noindent
We remember briefly how to construct a generic subordinated chain.\\

\noindent
A \emph{random point process} on the positive half-line is a sequence $\{T_n\}_{n\geq 0}$ of nonnegative random variable such that,
almost surely,

\begin{align*}
&T_0=0\\&
0<T_1<T_2<...\\&
\lim_{n\uparrow+\infty}T_n=+\infty
\end{align*}

\begin{definition}\label{def:unifmc}
Let $\{\hat{X}_n\}_{n\geq 0}$ be a discrete-time HMC with countable state space $E$ and transition matrix $\textbf{K}=\{k_{i,j}\}_{i,j\in E}$ and let
$\{T_n\}_{n\geq 1}$ be an HPP on $R_+$ with intensity $\lambda$ and associated counting process $N$. Suppose that $\{\hat{X}_n\}_{n\geq 0}$ and $\{T_n\}_{n\geq 1}$ are independent. The process $\{X(t)\}_{t\geq 0}$ with value in $E$ defined by

\begin{equation}\label{def:unifmc}
X(t)=\hat{X}_{N(t)}
\end{equation}

\noindent
is called \emph{uniform Markov chain}. The Poisson process $N$ is called the clock, and the chain $\{\hat{X}_n\}_{n\geq 0}$ is called the
\emph{subordinated chain}.
\end{definition}

\bigskip
\noindent
Let $\{X_n\}_{n\geq 0}$ a homogeneous discrete Markov chain with countable state space $N_0$, considering the transition probabilities

\begin{align}\label{equ:tmd}
\left\{
\begin{array}{ll}
P(X_{n+1}=j|X_n=i)=\frac{1}{1+i},\quad j=i+1\\
P(X_{n+1}=j|X_n=i)=\frac{i}{1+i},\quad j=i
\end{array}
\right.
\end{align}

\noindent
then from definition (\ref{def:unifmc}) $\{X_n\}_{n\geq 0}$ is the \emph{uniform Markov chain} of the birth process $\{X(t)\}_{t\geq 0}$ with associated differential equation (\ref{equ:bc}).\\

\noindent
Let $p_{n,k}=P(X_{n}=k)$ from (\ref{equ:tmd}) we obtain the following equations:

\begin{align}\label{equ:bd}
\left\{
\begin{array}{ll}
p_{n+1,k}=\frac{k}{k+1}p_{n,k}+\frac{1}{k}p_{n,k-1}, \quad \forall n\geq 0, \quad 1\leq k \leq n+1\\
p_{i,j}=0,\quad i\geq 0,\quad j\geq i+1\\
p_{0,0}=1,\quad p_{n,0}=0,\quad n\geq 1
\end{array}
\right.
\end{align}

\begin{proposition}\label{prop:aikfinalexpression}

\noindent
Let $p_{n,k}$ as in (\ref{equ:bd}) then

\begin{equation}\label{equ:ankfinal}
p_{n,k}=\frac{1}{k!}\sum_{i=1}^{k}(-1)^{k-i}\binom{k+1}{i+1}i^{k}\left(\frac{i}{i+1}\right)^{n-k}
\end{equation}

\end{proposition}

\noindent
\textbf{Proof}\\

\noindent
See Appendix.

\begin{flushright}
$\Box$
\end{flushright}

\subsection{Urn interpretation via dependent Bernoullian scheme}

\noindent
Consider an urn containing only \textbf{one} white ball and an arbitrary number of red balls.
If we draw the white ball then we add a red ball, while if we draw a red ball we do not anything.
In both cases we reinsert the drawn ball and proceed to the next drawing.\\
\noindent
This urn is a special type of Polya urn, see \cite{Janson}, corresponding to the following replacement scheme:

$$
\left(
\begin{array}{cc}
0 & 1 \\
0 & 0 \\
\end{array}\right)$$

\noindent
This Polya urn is assimilable to a sequence of dependent Bernoulli random variables, in fact if we consider that the number of red balls
in the urn corresponds to the number of previous successes then the probability of next success depends only on this number. 

\begin{proposition}\label{prop:urnscheme}

\noindent
Let $Z_n=\sum_{i=0}^{n}Y_i$, $n\geq 0$, where

\begin{equation}\label{berdepscheme}
Y_n|Y_0,...,Y_{n-1}\sim Ber\left(\frac{1}{1+\sum_{j=0}^{n-1}Y_j}\right);\quad Y_0=0
\end{equation}

\noindent
then $Z_n\stackrel{d}{=}X_n$, where the probability distribution of $X_n$ is given by (\ref{equ:tmd}).\\

\noindent
\textbf{Proof}\\

\end{proposition}

\noindent
It is sufficient to observe that the distribution of $Z_n$ satisfies (\ref{equ:bd}).

\begin{flushright}
$\Box$
\end{flushright}

\noindent
\textbf{Remark}\\

\noindent
The limit distribution of $Z_n$ is normal, see \cite{Janson} theorem 1.5.\\

\noindent
Now we recover the mean of the subordinated chain, i.e. the mean number of red balls in the urn after $n$ steps.\\
\noindent
Let $\{\tau_i\}_{i\geq 1}$ be a increasing sequence of positive r.v. such that $Y_n=1$ if and only if $n=\tau_i$ for some $i$.\\
\noindent
We have that:

\begin{equation*}
P(\tau_i=n)=P(Y_{\tau_{i-1}+1}=0,...,Y_{\tau_{i}-1}=0,Y_{n}=1)
\end{equation*}

\noindent
where $Y_m\sim Ber(\frac{1}{i})$, with $m\in[\tau_{i-1}+1,n]$.\\
\noindent
Then we have that $\tau_i\sim Geom(\frac{1}{i})$.\\
\noindent
Knowing that $P(X_n=k)=P(\sum_{i=1}^{n}Y_i=k)\Longleftrightarrow P(\sum_{i=1}^{k}\tau_i\leq n)$,
if we consider the r.v. $T=\sum_{i=1}^{k}\tau_i$, i.e. the number of trials necessary to obtain $k$ successes, we have

\begin{equation*}
E(T)=E\left(\sum_{i=1}^{k}\tau_i\right)=\sum_{i=1}^{k}i
\end{equation*}

\noindent
Solving the equation $E(T)=n$, with respect to the variable $k$, i.e. $\frac{k(k+1)}{2}=n$, we find

\begin{equation}\label{equ:dmean}
E(X_n)=\frac{-1+\sqrt{1+8n}}{2}
\end{equation}

\noindent
We want to give a result of weak concentration for the embedded process $X_n$. For this purpose we give an upper bound for the variance.\\
\noindent
Consider (\ref{berdepscheme}), for the second moment of $X_n$ holds

\begin{align*}
E(X_{n+1}^2)&=E(E(X_{n+1}^2|X_n))=E(E((X_{n}+Y_{n+1})^2|X_n))=\\&
=E\left(X^2_n + \frac{1}{1+X_n} +2X_n\frac{1}{1+X_n}\right)=\\&
=1+E(X^2_n)+E\left(\frac{X_n}{1+X_n}\right)
\end{align*}

\noindent
so that

\begin{align*}
E(X^2_{n+1})\leq 2 +E(X^2_n)
\end{align*}

\noindent
Being $E(X^2_1)\leq 1$ we obtain $E(X^2_{n})\leq 2n$ then

\begin{equation}\label{equ:var_ank}
Var(X_{n})\leq 2n-\left(\frac{-1+\sqrt{1+8n}}{2}\right)^2=-\frac{1}{2}+\sqrt{1+8n}
\end{equation}

\begin{proposition}
\noindent
Let $\{X_n\}_{n\geq 0}$ the process with probability distribution function as in (\ref{equ:bd}), then for $n\rightarrow +\infty$

\begin{equation}\label{equ:cheb_ank}
P(|X_n-E(X_n)|\geq \epsilon E(X_n))\leq \frac{Var(X_n)}{\epsilon^2 E^2(X_n)}
\sim \frac{\sqrt{2}}{\epsilon^2 \sqrt{n}}\rightarrow 0
\end{equation}

\noindent
\textbf{Proof}\\

\end{proposition}

\noindent
It is an obvious consequence of (\ref{equ:dmean}), (\ref{equ:var_ank}) and Chebychev's inequality.

\begin{flushright}
$\Box$
\end{flushright}

\noindent
Considering the Bernoullian scheme representation (\ref{prop:urnscheme}) for the subordinated chain, we establish the large deviation bounds
using some results due to McDiarmid \cite{McDiarmidcentering}.\\
\noindent
To this purpose we recall the definition of centering sequences and some result about them.\\
\noindent
For more details see \cite{McDiarmidcentering}.\\

\bigskip

\begin{definition}\label{def:centering}
\noindent
Given a sequence $X =(X_1,X_2,...)$ of (integrable) random variables the corresponding difference sequence is $Y =(Y_1,Y_2,...)$ where $Y_k =X_k-X_{k-1}$ (and where we always set $X_0=0$). Let $\mu_k(x) = E(Y_k|X_{k-1} = x)$, that is $\mu_k(X_{k-1})$ is a version of $E(Y_k|X_{k-1})$.
We call the sequence $X$ centering if for each $k =2,3,...$ we may take $\mu_k(x)$ to be a non-increasing function of $x$.
\end{definition}

\begin{theorem}\label{thm:centeringdeviation}
\noindent
Let $X_1,X_2,...,X_n$ be a centering sequence with\\
\noindent
corresponding differences $Y_k = X_k-X_{k-1}$ satisfying $0\leq Y_k\leq 1$ for
each k. Then

\begin{align*}
&P(X_n\geq(1+\epsilon)E(X_n))\leq \exp\left(-\frac{1}{3}\epsilon^2 E(X_n)\right) \quad 0<\epsilon <1;\\&
P(X_n\leq(1-\epsilon)E(X_n))\leq \exp\left(-\frac{1}{3}\epsilon^2 E(X_n)\right) \quad 0< \epsilon <1
\end{align*}

\noindent
\textbf{Proof}\\

\end{theorem}

\noindent
See \cite{McDiarmidcentering}.

\begin{flushright}
$\Box$
\end{flushright}

\noindent
\textbf{Remark}\\

\noindent
Considering the Bernoullian scheme as in proposition (\ref{prop:urnscheme}),
let $Y_n=X_n-X_{n-1}$ we have $E(Y_n|X_{n-1})=\frac{1}{1+X_{n-1}}$, for all $n\geq 1$, then $X_n$ is a centering sequence.\\

\begin{corollary}
Under the hypothesis of proposition (\ref{prop:urnscheme}) we have

\begin{align}\label{largedeviationb}
&P(X_n\geq(1+\epsilon)E(X_n))\leq \exp\left(-\frac{1}{3}\epsilon^2 E(X_n)\right)\sim e^{-\frac{1}{3}\epsilon^2\sqrt{2n}}, \quad 0<\epsilon <1\\&
P(X_n\leq(1-\epsilon)E(X_n))\leq \exp\left(-\frac{1}{3}\epsilon^2 E(X_n)\right)\sim e^{-\frac{1}{3}\epsilon^2\sqrt{2n}}, \quad 0< \epsilon <1
\end{align}

\noindent
\textbf{Proof}\\

\end{corollary}

\noindent
The thesis is an obvious consequence of (\ref{equ:dmean}) and theorem  (\ref{thm:centeringdeviation}).

\begin{flushright}
$\Box$
\end{flushright}

\noindent
\begin{proposition}
\noindent
Let $X_1,X_2,...,X_n$ be a centering sequence with corresponding differences $Y_k = X_k-X_{k-1}$ with mean $\mu_k=Y_k$ and suppose that there are
constants $a_k$ and $b_k$ such that $a_k\leq Y_k \leq b_k$ for each $k$. Then for any $h > 0$

\begin{align}\label{mgfbound}
E(exp(hX_n))\leq \prod_{k=1}^{n} \left(\frac{b_k-\mu_k}{b_k-a_k}e^{h a_k}+\frac{\mu_k-a_k}{b_k-a_k}e^{h b_k}\right)
\end{align}

\noindent
\textbf{Proof}\\

\end{proposition}

\noindent
See \cite{McDiarmidcentering}.

\begin{flushright}
$\Box$
\end{flushright}

\begin{corollary}
Under the hypothesis of above proposition, if $a_k=0$ and $b_k=1$ then

\begin{align}\label{mgfbound}
E(exp(hX_n))\leq (1-\alpha+\alpha e^h)^n
\end{align}

\noindent
where $\alpha=\frac{E(X_n)}{n}$.\\

\noindent
\textbf{Proof}\\

\end{corollary}

\noindent
See \cite{McDiarmidcentering}.

\begin{flushright}
$\Box$
\end{flushright}

\noindent
\textbf{Remark}\\

\noindent
From (\ref{equ:dmean}) we have $\alpha\sim \frac{\sqrt{2}}{\sqrt{n}}$ recalling (\ref{mgfbound}) then for $n\rightarrow +\infty$

\begin{align}\label{mgfbpbound}
E(exp(hX_n))\sim e^{\sqrt{2}(e^h-1)\sqrt{n}}
\end{align}

\newpage
\section{Appendix}

\noindent
First we rewrite (\ref{equ:bd}) as

\begin{equation}\label{ank_hsb}
p_{n,k}=p_{n,k-1}+\left(k+\frac{1}{k+1}\right)p_{n,k}-kp_{n+1,k+1}, \quad \forall n\geq 0, \quad 0\leq k \leq n+1
\end{equation}

\noindent
where

\begin{equation*}
p_{n,-1}=0,\quad \forall n\geq 0
\end{equation*}

\noindent
We define

\begin{equation}\label{equ:fgd}
F(x,y)=\sum_{n=0}^{+\infty}\sum_{k=0}^{+\infty}x^n y^kp_{n,k}=\sum_{n=0}^{+\infty}\sum_{k=0}^{n}x^n y^kp_{n,k}
\end{equation}

\noindent
substituting into (\ref{ank_hsb}):

\begin{align}\label{equ:bgenfun}
\notag\sum_{n=0}^{+\infty}\sum_{k=0}^{n}x^n y^kp_{n,k}=&\sum_{n=0}^{+\infty}\sum_{k=0}^{n}x^n y^kp_{n,k-1}+\sum_{n=0}^{+\infty}\sum_{k=0}^{n}x^n y^kkp_{n,k}+\\&
-\sum_{n=0}^{+\infty}\sum_{k=0}^{n}x^n y^kkp_{n+1,k}+\sum_{n=0}^{+\infty}\sum_{k=0}^{n}x^n y^k\frac{1}{1+k}p_{n,k}
\end{align}

\noindent
After the computation of the all summands of (\ref{equ:bgenfun}), we recover that:\\

\begin{align}\label{diffgdfirst}
F(x,y)=yF(x,y)+yF_y(x,y)-\frac{y}{x}F_y(x,y)+\frac{1}{y}\int_{0}^{y}F(x,t)dt
\end{align}

\noindent
Now from (\ref{diffgdfirst}) multiplying both members by $y$ and computing the derivative respect to $y$ we have

\begin{equation}\label{equ:difffgd}
y \frac{1-x}{x} F_{yy}(x,y)+\left(\frac{2-x}{x}-y\right) F_y(x,y)-2F(x,y)=0
\end{equation}

\begin{proposition}\label{equ:diffgenfuncbp}

\noindent
For $0\leq x <1$, $0\leq y < 1$ let

\begin{equation}\label{equ:finalF}
F(x,y)=1+\sum_{k=1}^{+\infty}\frac{x^ky^k}{k!\prod_{h=0}^{k}(1-\frac{hx}{h+1})}
\end{equation}

\noindent
then (\ref{equ:finalF}) is a solution of (\ref{equ:difffgd}) verifying the conditions $F(0,y)=F(x,0)=1$.\\

\noindent
\textbf{Proof}\\

\end{proposition}

\noindent
From (\ref{equ:difffgd}) we note that the variable $x$ does not appear in the derivatives, therefore we treat it as a parameter.\\
\noindent
For convenience we put, for $0<x<1$,

\begin{align}\label{abconstantbp}
&a=\frac{1-x}{x}; \quad b=\frac{2-x}{x}
\end{align}

\noindent
then (\ref{equ:difffgd}) becomes

\begin{equation}\label{equ:ndifffgd}
F_{yy}(x,y)+\frac{b-y}{ay}F_y(x,y)-\frac{2}{ay}F(x,y)=0
\end{equation}

\noindent
We follow Frobenius method. The solutions of the indicial equation

\begin{equation*}
\lambda(\lambda-1)+\lambda\frac{b}{a}=0
\end{equation*}

\bigskip
\noindent
are $\lambda_1=0$ and $\lambda_2=1-\frac{b}{a}$.\\
\noindent
Since $\lambda_2<\lambda_1$ then all the solutions of (\ref{equ:ndifffgd}) will be of the type

\begin{equation*}
F(x,y)=C_1(x)F_1(y)+C_2(x)F_2(y)
\end{equation*}

\noindent
where

\begin{align*}
&F_1(y)=\sum_{k=0}^{+\infty}c_ky^k,\quad c_0=1\\&
F_2(y)=y^{1-\frac{b}{a}}\sum_{k=0}^{+\infty}d_ky^k+ C F_1(y)\ln(y),\quad d_0=1
\end{align*}

\noindent
where the constant $C$ is equal to zero if $\lambda_1-\lambda_2$ is not an integer, i.e. $x\neq 1-\frac{1}{n}$.\\
\noindent
After computing the derivatives of $F_1(y)$ and substituting in (\ref{equ:difffgd}) we recover that

\begin{align*}
\left\{
\begin{array}{ll}
c_0=1\\
c_k=\frac{k+1}{\prod_{i=0}^{k-1}(ia+b)},\quad \forall k\geq 1
\end{array}
\right.
\end{align*}

\noindent
thus

\begin{align*}
F_1(y)=1+\sum_{k=1}^{+\infty}\frac{k+1}{\prod_{i=0}^{k-1}(ia+b)}y^k
\end{align*}

\noindent
From (\ref{abconstantbp}) we have

\begin{equation*}
F(x,y)\doteq 1+\sum_{k=1}^{+\infty}\frac{x^ky^k}{k!\prod_{h=0}^{k}(1-\frac{hx}{h+1})}
\end{equation*}

\noindent
is well defined even if $x=0$ and we have $F(0,y)=F(x,0)=1$.\\
\noindent
This conclude the proof.

\begin{flushright}
$\Box$
\end{flushright}

\begin{proposition}\label{equ:pnkexplicit}

\noindent
Let $F(x,y)$ as in (\ref{equ:finalF}) then

\begin{align}\label{equ:newfinalF}
&F(x,y)=1+\sum_{n=1}^{+\infty}\sum_{k=1}^{n}x^{n}y^k\frac{1}{k!}\sum_{i=1}^{k}A_{i,k}\left(\frac{i}{i+1}\right)^{n-k}
\end{align}

\noindent
where

\begin{align}\label{equ:aikfinalF}
A_{i,k}=\frac{1}{\prod_{h=1,h\neq i}^{k}(1-\frac{h}{h+1}\frac{i+1}{i})}
\end{align}

\noindent
\textbf{Proof}\\

\end{proposition}

\noindent
Let

\begin{align*}
\left\{
\begin{array}{ll}
\phi_1=\frac{1}{1-\frac{x}{2}}\\
\phi_k=\frac{1}{\prod_{h=0}^{k}\left(1-\frac{hx}{h+1}\right)}, \quad k>1
\end{array}
\right.
\end{align*}

\noindent
because

\begin{equation*}
\phi_k=\sum_{i=1}^{k}\frac{A_{i,k}}{1-\frac{ix}{i+1}}
\end{equation*}

\noindent
then

\begin{equation*}
\sum_{i=1}^{k}A_{i,k}\prod_{h=1,h\neq i}^{k}(1-\frac{hx}{h+1})=1
\end{equation*}

\noindent
Choosing $x=\frac{i}{i+1}$ we obtain

\begin{equation*}
A_{i,k}=\frac{1}{\prod_{h=1,h\neq i}^{k}(1-\frac{h}{h+1}\frac{i+1}{i})}
\end{equation*}

\noindent
From (\ref{equ:newfinalF})

\begin{align*}
F(x,y)&=1+\sum_{k=1}^{+\infty}x^ky^k\frac{1}{k!}\sum_{i=1}^{k}\frac{A_{i,k}}{1-\frac{ix}{i+1}}=\\&
=1+\sum_{k=1}^{+\infty}x^ky^k\frac{1}{k!}\sum_{i=1}^{k}A_{i,k}\sum_{r=0}^{+\infty}\left(\frac{i}{i+1}\right)^rx^r=
\\&=1+\sum_{k=1}^{+\infty}x^ky^k\frac{1}{k!}\sum_{i=1}^{k}A_{i,k}\sum_{s=1}^{+\infty}\left(\frac{i}{i+1}\right)^{s-1}x^{s-1}
\end{align*}

\noindent
Inverting the order of the last summations it follows that

\begin{align*}
F(x,y)&=1+\sum_{k=1}^{+\infty}x^ky^k\frac{1}{k!}\sum_{s=1}^{+\infty}x^{s-1}\sum_{i=1}^{k}A_{i,k}\left(\frac{i}{i+1}\right)^{s-1}=\\&
=1+\sum_{k=1}^{+\infty}\sum_{s=1}^{+\infty}x^{k+s-1}y^k\frac{1}{k!}\sum_{i=1}^{k}A_{i,k}\left(\frac{i}{i+1}\right)^{s-1}=
\\&=1+\sum_{k=1}^{+\infty}\sum_{n=k}^{+\infty}x^{n}y^k\frac{1}{k!}\sum_{i=1}^{k}A_{i,k}\left(\frac{i}{i+1}\right)^{n-k}
\end{align*}

\noindent
Inverting the order of the first two summations we have

\begin{align*}
&F(x,y)=1+\sum_{n=1}^{+\infty}\sum_{k=1}^{n}x^{n}y^k\frac{1}{k!}\sum_{i=1}^{k}A_{i,k}\left(\frac{i}{i+1}\right)^{n-k}
\end{align*}

\noindent
The thesis follows.

\begin{flushright}
$\Box$
\end{flushright}

\noindent
\textbf{Remark}\\

\noindent
From (\ref{equ:newfinalF}) taking into account (\ref{equ:fgd}) for $n\geq k\geq 1$ it follows that

\begin{equation}\label{equ:ank1}
p_{n,k}=\frac{1}{k!}\sum_{i=1}^{k}A_{i,k}\left(\frac{i}{i+1}\right)^{n-k}
\end{equation}

\noindent
Now we conclude proving that (\ref{equ:ank1}) is the solution of the initial system (\ref{equ:bd}).

\begin{proposition}\label{prop:pnkverify}

\noindent
Let $p_{n,k}$ as in (\ref{equ:ank1}) then it solves (\ref{equ:bd}).\\

\noindent
\textbf{Proof}\\

\end{proposition}

\noindent
We proceed by induction. Let $k=1$ in this case we have

\begin{equation*}
p_{n+1,1}=\frac{1}{2}p_{n,1}
\end{equation*}

\noindent
from above equality substituting in (\ref{equ:ank1}) it follows that

\begin{align*}
&\frac{1}{1!}\sum_{i=1}^{1}A_{i,1}\left(\frac{i}{i+1}\right)^{n+1-1}=\frac{1}{2}\frac{1}{1!}\sum_{i=1}^{1}A_{i,1}\left(\frac{i}{i+1}\right)^{n-1}\\&
\mbox{i.e.}\\&
A_{1,1}\left(\frac{1}{2}\right)^n=\frac{1}{2}A_{1,1}\left(\frac{1}{2}\right)^{n-1}
\end{align*}

\noindent
then the thesis, remembering that $A_{1,1}=1$.\\

\noindent
For $k=2,...,n$ and $n\geq 1$ substituting in (\ref{equ:ank1}) we have

\begin{equation*}
\frac{1}{k!}\sum_{i=1}^{k}A_{i,k}\left(\frac{i}{i+1}\right)^{n+1-k}=\frac{k}{k+1}\frac{1}{k!}\sum_{i=1}^{k}A_{i,k}\left(\frac{i}{i+1}\right)^{n-k}+
\frac{1}{k}\frac{1}{(k-1)!}\sum_{i=1}^{k-1}A_{i,k-1}\left(\frac{i}{i+1}\right)^{n-k+1}
\end{equation*}

\noindent
multiplying both members for $k!$ and eliminate by algebraic calculation the terms with index $i=k$ we have

\begin{align*}
\sum_{i=1}^{k-1}A_{i,k}\left(\frac{i}{i+1}\right)^{n+1-k}=\frac{k}{k+1}\sum_{i=1}^{k-1}A_{i,k}\left(\frac{i}{i+1}\right)^{n-k}+
\sum_{i=1}^{k-1}A_{i,k-1}\left(\frac{i}{i+1}\right)^{n-k+1}
\end{align*}

\noindent
So that we have to verify that

\begin{align}\label{equ:aikpnkverify}
\sum_{i=1}^{k-1}\left(\frac{i}{i+1}\right)^{n-k}\left[A_{i,k}\left(\frac{i}{i+1}\right)-\frac{k}{k+1}A_{i,k}-A_{i,k-1}\left(\frac{i}{i+1}\right)\right]=0
\end{align}

\noindent
From (\ref{equ:aikpnkverify}) for $k=2$ we have

\begin{align*}
&\sum_{i=1}^{1}\left(\frac{i}{i+1}\right)^{n-2}\left[A_{i,2}\left(\frac{i}{i+1}\right)-\frac{2}{2+1}A_{i,2}-A_{i,1}\left(\frac{i}{i+1}\right)\right]=0
\end{align*}

\noindent
because

\begin{align*}
&A_{1,2}=\frac{1}{\prod_{h=1,h\neq i}^{2}(1-\frac{h}{h+1}2)}=-3
\end{align*}

\noindent
and remembering that $A_{1,1}=1$.\\

\noindent
For $k>2$, let us observe that for every $i=1,...,k-1$ we have

\begin{align*}
&\left(\frac{i}{i+1}-\frac{k}{k+1}\right)\frac{1}{\prod_{h=1,h\neq i}^{k}\left(1-\frac{h}{h+1}\frac{i+1}{i}\right)}
-\left(\frac{i}{i+1}\right)\frac{1}{\prod_{h=1,h\neq i}^{k-1}(1-\frac{h}{h+1}\frac{i+1}{i})}=\\& \phantom{}\\&
=\left(\frac{i}{i+1}-\frac{k}{k+1}\right)\frac{1}{[\prod_{h=1,h\neq i}^{k-1}(1-\frac{h}{h+1}\frac{i+1}{i})]}
\frac{1}{(1-\frac{k}{k+1}\frac{i+1}{i})}
-\left(\frac{i}{i+1}\right)\frac{1}{\prod_{h=1,h\neq i}^{k-1}\left(1-\frac{h}{h+1}\frac{i+1}{i}\right)}=\\&
\phantom{}\\&
=\frac{1}{\prod_{h=1,h\neq i}^{k-1}(1-\frac{h}{h+1}\frac{i+1}{i})}\left[\frac{i}{i+1}-\frac{k}{k+1}-\frac{i}{i+1}\left(1-\frac{k}{k+1}\frac{i+1}{i}\right)\right]=\\&
\phantom{}\\&
=\frac{1}{\prod_{h=1,h\neq i}^{k-1}(1-\frac{h}{h+1}\frac{i+1}{i})}\left[\frac{i}{i+1}-\frac{k}{k+1}-\frac{i}{i+1}+\frac{k}{k+1}\right]=0
\end{align*}

\noindent
then the thesis follows.

\begin{flushright}
$\Box$
\end{flushright}

\noindent
We proceed to establish a compact form for the $p_{n,k}$.\\

\begin{proposition}\label{prop:aikfinalexpression}

\noindent
Under the hypothesis of proposition (\ref{prop:pnkverify}) for $n\geq k\geq1$

\begin{equation}\label{equ:ankfinal}
p_{n,k}=\frac{1}{k!}\sum_{i=1}^{k}(-1)^{k-i}\binom{k+1}{i+1}i^{k}\left(\frac{i}{i+1}\right)^{n-k}
\end{equation}

\noindent
\textbf{Proof}\\

\end{proposition}

\noindent
From (\ref{equ:ank1}) we have that

\begin{align*}
\left\{
\begin{array}{ll}
A_{1,1}=1\\
A_{i,k}=\frac{1}{\prod_{h=1,h\neq i}^{k}(1-\frac{h}{h+1}\frac{i+1}{i})}, \quad \forall k>1
\end{array}
\right.
\end{align*}

\noindent
so that

\begin{align*}
\frac{1}{A_{i,k}}&=\prod_{h=1,h\neq i}^{k}\left(1-\frac{h}{h+1}\frac{i+1}{i}\right)=\\&
=\left(\frac{i+1}{i}\right)^{k-1}\prod_{h=1,h\neq i}^{k}\frac{i-h}{(i+1)(h+1)}=\\&
=\frac{1}{i^{k-1}}\frac{i+1}{(k+1)!}\prod_{h=1,h\neq i}^{k}(i-h)
\end{align*}

\noindent
Let $i=k$: then

\begin{align*}
&\frac{1}{A_{k,k}}=\frac{1}{k^{k-1}}\frac{k+1}{(k+1)!}\prod_{h=1,h\neq k}^{k}(k-h)=
\frac{1}{k^{k-1}}\frac{k+1}{(k+1)!}(k-1)!=\frac{1}{k^{k}}
\end{align*}

\noindent
Let $i<k$: then

\begin{align*}
\frac{1}{A_{i,k}}&=\frac{1}{i^{k-1}}\frac{i+1}{(k+1)!}\prod_{h=1,h\neq i}^{k}(i-h)=\\&
=\frac{1}{i^{k-1}}\frac{i+1}{(k+1)!}\prod_{h=1}^{i-1}(i-h)\prod_{h=i+1}^{k}(i-h)=
\frac{1}{i^{k-1}}\frac{i+1}{(k+1)!}(i-1)!\prod_{h=i+1}^{k}(i-h)=\\&
=\frac{1}{i^{k-1}}\frac{i+1}{(k+1)!}(i-1)!(-1)^{k-i}\prod_{h=i+1}^{k}(h-i)=\frac{1}{i^{k}}\frac{(i+1)!}{(k+1)!}(-1)^{k-i}\prod_{h=i+1}^{k}(h-i)=
\\&=\frac{1}{i^{k}}\frac{(i+1)!}{(k+1)!}(-1)^{k-i}\prod_{h-i=1}^{k}(h-i)=\frac{1}{i^{k}}\frac{(i+1)!}{(k+1)!}(-1)^{k-i}(k-i)!=\\&
=\frac{1}{i^{k}}(-1)^{k-i}\frac{1}{\binom{k+1}{i+1}}
\end{align*}

\noindent
so that

\begin{align*}
A_{i,k}=(-1)^{k-i}\binom{k+1}{i+1}i^{k},\quad \forall k>1,\quad i=1,...,k
\end{align*}

\noindent
then the thesis follows.

\begin{flushright}
$\Box$
\end{flushright}


\begin{thebibliography}{10}
\bibitem[1]{Fellervol1}Feller, W., An introduction to probability theory and its applications. {V}ol. {I}, Third edition, John Wiley \& Sons Inc., New York, 1968
\bibitem[2]{McDiarmidcentering} McDiarmid, C., Centering sequences with bounded differences, Combin. Probab. Comput., Vol.6, 1997, No.1, 79--86
\bibitem[3]{Janson} Janson, S., Limit theorems for triangular urn schemes, Probab. Theory Related Fields, Vol. 134, 2006, No.3, pp. 417-452.
\end{thebibliography}
\end{document}